\documentclass[11pt,a4paper]{amsart}
\usepackage[english]{babel}
\usepackage[utf8]{inputenc} 
\usepackage[T1]{fontenc}
\usepackage{verbatim}
\usepackage{palatino}
\usepackage{amsmath}
\usepackage{mathabx}
\usepackage{amsbsy}
\usepackage{amssymb}
\usepackage{amsthm}
\usepackage{amsfonts}
\usepackage{graphicx}
\usepackage{esint}
\usepackage{color}
\usepackage{mathtools}
\usepackage{overpic}
\usepackage{mathrsfs}
\usepackage{tabularx}
\usepackage{datetime}
\renewcommand{\today}{\the\day\ \shortmonthname[\month] \the\year}
\subjclass[2010]{05B99, 28A78, 28A80}
\keywords{Discretised sum-product, discretised ring theorem, Shannon entropy, Ahlfors-regular sets}
\thanks{The author is supported in part by an NSERC Alliance grant administered by Pablo Shmerkin and Joshua Zahl.}
\numberwithin{equation}{section}

\theoremstyle{plain}

\newtheorem{theorem}[equation]{Theorem}
\newtheorem{proposition}[equation]{Proposition}

\newtheorem{lemma}[equation]{Lemma}

\newtheorem{corollary}[equation]{Corollary}

\theoremstyle{definition}

\newtheorem{definition}[equation]{Definition}
\newtheorem{remark}[equation]{Remark}

\mathtoolsset{showonlyrefs}

\newcommand{\N}{\mathbb{N}}
\newcommand{\Z}{\mathbb{Z}}
\newcommand{\pr}{\mathbb{P}}

\newcommand{\rmC}{\mathrm{C}}
\newcommand{\calF}{\mathcal{F}}
\newcommand{\diam}{\operatorname{diam}}

\newcommand{\D}{\mathcal{D}}
\newcommand{\R}{\mathbb{R}}
\newcommand{\E}{\mathrm{H}}
\newcommand{\Ed}{\mathrm{H}_\delta}
\newcommand{\col}{\operatorname{col}_\delta}
\newcommand{\collis}{\operatorname{col}}

\newcommand{\spt}{\operatorname{spt}}

\title{Sum-product phenomena for  Ahlfors-regular sets}
\author{William O'Regan}
\date{\today}
\address{Department of Mathematics, The University of British Columbia, Room 121, 1984 Mathematics Road, Vancouver  BC Canada V6T 1Z2}
\email{woregan@math.ubc.ca}
\begin{document}
\maketitle
\begin{abstract}
    We utilise the recent work of Orponen to yield a sum-product result for Ahlfors-regular sets. As a corollary, we obtain the fractal analogue of Solymosi's $4/3$-bound for finite subsets of $\R.$
\end{abstract}
\section{Introduction}
The sum-product phenomenon is the maxim that additive and multiplicative structure find it hard to co-exist. An example of this is the result of Erd\H{o}s and Szemer\'edi \cite{erdsz}: There exists $\epsilon >0$ so that for all large enough $A \subset \Z$ we have
\begin{equation}
    |A+A| + |AA| > |A|^{1+\epsilon}.
\end{equation}
In other words, finite subsets of $\Z$ cannot closely resemble rings. The above is now known to hold over $\R$ and the best known value of $\epsilon$ (in print) one can take is due to Rudnev and Stevens \cite{rudnev2022update}; $\epsilon = 1/3 + 2/1167.$ (See also the recent preprint \cite{bloom}  that gives a modest improvement to the work of Rudnev--Stevens). It is conjectured that one may take any $0< \epsilon <1,$ in other words, either the sum-set or the product-set must be almost as large as possible. 

\indent Sum-product phenomena are now known to hold in a large variety of settings, in particular, for fractal sets. Loosely stated, if $A \subset \R$ has `dimension' $s$, then one of $A+A$ or $AA$ has `dimension' `much larger' than $s.$ This problem was introduced by Katz and Tao \cite{katz2001some} as it is related to the now solved Erd\H{o}s--Volkmann ring problem \cite{erd}: are there Borel subrings of $\R$ with Hausdorff dimension strictly between $0$ and $1?.$ This was answered in the negative by Bourgain \cite{bou03} and Edgar--Miller \cite{edgar2003borel} independently. 

\indent The Edgar--Miller paper gave a very direct and fairly elementary proof. Bourgain showed that the ring problem is a corollary of the so called `discretised ring theorem': For all $0 < s < 1$ there exists $\epsilon > 0$ so that if $A \subset \R$ is a finite $\delta$-separated set that resembles a fractal set \footnote{$|A| \approx \delta^{-s}, |A \cap B(x,r)| \lessapprox r^s|A|$ for all $r \geq \delta, x\in \R$ }of dimension $s,$ then
\begin{equation}
    N_\delta(A+A) + N_\delta(AA) > |A|^{1+\epsilon}
\end{equation}
provided that $\delta > 0$ is small enough. Here, and throughout, $N_\delta$ denotes the least number of closed intervals of length $\delta$ needed to cover the set. This problem has attracted a large amount of interest in recent years. For the best bound on $\epsilon$ see \cite{furen}, \cite{renwang}. For an elementary proof see \cite{gut}. For a result with weaker conditions on $A$ see \cite{bougam}, \cite{bou}. 

\indent The aim of this note is to improve on the bound given in \cite{renwang}, but for the more restricted class of regular sets. See the below definition.
\begin{definition}
    Let $C,s > 0.$ A set $K \subset \R^d$ is called \textit{upper $(s,C)$-regular} if 
    \begin{equation}
        N_r(K \cap B(x,R)) \leq C\big(\tfrac{R}{r}\big)^s
    \end{equation}
    for all $0 < r \leq R < \infty$ and $x \in \R^d.$ 
    
    \indent A Borel probability measure $\mu$ is $(s,C)$-\textit{Frostman} if $\mu(B(x,r)) \leq Cr^s$ for all $r >0$ and $x \in \R^d.$ Further, it is $(s,C)$-\textit{regular} if it is $(s,C)$-Frostman and $\spt \mu$ is upper $(s,C)$-regular. 
    
    \indent A set $K \subset \R^d$ is \textit{$(s,C)$-regular} if there exists an $(s,C)$-regular measure $\mu$ with $K = \spt \mu.$
\end{definition}
\begin{theorem}\label{thm.main}
    Let $0 < s \leq 1/2, \eta >0$ and let $C > 0.$ There exists a $\delta_0 = \delta_0(C,s,\eta) >0$ so that the following holds. 
Let $\mu$ be an $(s,C)$-regular measure with $\spt \mu \subset [1,2],$ and let $X,Y$ be iid random variables distributed by $\mu.$ We have
\begin{equation}
    \Ed(X+Y) + 2\Ed(XY) > (4s -\eta)\log(1/\delta) - O_C(1),
\end{equation}
    for all $0 < \delta < \delta_0.$ 
    \end{theorem}
     Here, and throughout, $\Ed$ denotes the Shannon entropy with respect to a partition of intervals of length $\delta.$ This will be defined and made precise in section \ref{section.prelim}. Theorem \ref{thm.main} is sharp in the sense that $4s$ cannot be replaced by a larger number; see section \ref{section.sharp}. We remark that in Theorem \ref{thm.main}, we may replace $X+Y$ with $X-Y$ and/or replace $XY$ with $X/Y,$ by minor alterations in the proof. An immediate corollary follows using the fact that $\Ed(\cdot) \leq \log N_\delta(\cdot) + O(1).$
\begin{theorem}\label{thm.cor}
    Let $0 < s \leq 1/2,\eta > 0$ and let $C > 0.$ There exists a $\delta_0 = \delta_0(C,s,\eta) >0$ so that the following holds. 
Let $A \subset [1,2]$ be an $(s,C)$-regular set. We have
\begin{equation}\label{eq.sharp}
    N_\delta(A+A)N_\delta(AA)^2 > \delta^{-4s + \eta},
\end{equation}
and
\begin{equation}\label{eq.soly}
    N_\delta(A+A) + N_\delta(AA) > \delta^{-4s/3+\eta}
\end{equation}
    for all $0 < \delta < \delta_0.$ 
    \end{theorem}
    Again, \eqref{eq.sharp} is sharp as we will see in section \ref{section.sharp}, and we may replace $A+A$ with $A-A$ and/ or replace $AA$ with $A/A.$ Inequality \eqref{eq.soly} resembles the bound obtained by Solymosi in \cite{sol} for the discrete sum-product problem, although the proofs seem totally unrelated. 
\subsection{Proof sketch}
Fix $0 < s \leq 1/2.$ Here we use $\gtrapprox$ loosely to morally mean greater than.  A remarkable recent result of Orponen \cite{orpahlfors} shows that for an $s$-regular measure $\mu$ and a $\epsilon$-Frostman measure $\nu,$ both supported on $\R,$ we may always find $x \in \spt \nu$ so that $$N_\delta(\spt \mu + x \spt \mu) \gtrapprox \delta^{-2s}.$$
As an artifact of his proof, we are able to massage and manipulate his work into the following result. Let $\mu$ be Ahlfors $s$-regular. Let $X,Y,Z$ be iidrandom variables distributed by $\mu.$ Then
\begin{equation}
    \Ed((X+Y)Z) \gtrapprox 2s\log(1/\delta).
\end{equation}
We may then use the submodularity of Shannon entropy and the fact that $X,Y,Z$ are iidto show that
\begin{equation}
    \Ed((X+Y)Z) \leq \Ed(X+Y) + 2\Ed(XY) - 2\Ed(X) + O(1).
\end{equation}
 Combining these inequalities gives the required result. The use of the submodularity of entropy has been previously utilised in the work of M\'ath\'e--O'Regan to gain a strong bound for the discretised ring theorem. See \cite{mator}.
\begin{remark}
The result we use \cite[Theorem 1.13]{orpahlfors}, is only valid for regular measures of dimension $\leq 1.$ If one updates this result appropriately to dimensions $ > 1,$ then Theorem \ref{thm.main} and Theorem \ref{thm.cor} may be suitably updated to include results for the full range $0 < s < 1.$ When $1/2 < s < 1$ the '$4s$' in Theorem \ref{thm.main} and Theorem \ref{thm.cor} would be replaced by $1 + 2s.$ See Remark \ref{rem.remark} for the version of \cite[Theorem 1.13]{orpahlfors} that would be required here. 
\end{remark}
\subsection*{Acknowledgements}
We thank Pablo Shmerkin and Joshua Zahl for invaluable discussions and suggestions, in particular to Joshua for pointing out that Theorem \ref{thm.main} is sharp. We thank Tuomas Orponen for useful comments. We thank the anonymous referee for comments improving the accuracy and quality of the exposition. 
\section{Preliminaries}\label{section.prelim}
\subsection{Entropy}
In the below our random variables may take finitely many values in $\R^d.$ Let $X$ be a random variable. Define the \textit{Shannon entropy} of $X$ by
\begin{equation}
    \E(X) = -\sum_{x}\pr(X=x)\log \pr(X=x).
\end{equation}
Define the \textit{collision entropy} of $X$ by
\begin{equation}
    \collis(X) = -\log \sum_{x}\pr(X=x)^2.
\end{equation} Here, and throughout, $\log$ will be taken to base 2, but this is not particularly important. We adhere to the convention that $0\log 0 = 0.$ 

\indent We recall some useful facts. The first is monotonicity:
For any random variable $X$ we have
\begin{equation}
    \collis(X) \leq \E(X).
\end{equation}
The second is that Shannon entropy is concave. The third is the chain rule. Let $X,Y$ be random variables, then
\begin{equation}\label{eq.chainrule}
    \E(X,Y)=\E(X)+\E(Y|X).
\end{equation}The below well-known inequality is useful. For example, see \cite[Lemma A.2]{taoent}. For two random variables $X,Y$ we say that $X$ \textit{determines} $Y$ if $\sigma(Y) \subset \sigma(X),$ where $\sigma(\cdot)$ is the $\sigma$-algebra generated by $\cdot.$
\begin{theorem}[Submodular inequality]\label{thm.submod}
    Let $X,Y,Z,W$ be random variables. Suppose that $Z$ determines $X$ and $W$ determines $X.$ Suppose that $(Z,W)$ determines $Y.$ Then
    \begin{equation}
        \E(X) + \E(Y) \leq \E(Z) + \E(W).
    \end{equation}
\end{theorem}
\indent We wish to not restrict ourselves to just random variables which have finite support. But, we do want to use the theory above. We do this by discretising our (infinitely supported) random variables at a scale $\delta >0.$ To this end, let $\D_\delta = \D_\delta(\R^d)$ be the collection of 
$\delta$-cubes of the form $[\delta i_1, \delta (i_1+1)) \times \cdots \times [\delta i_d, \delta (i_d +1 )), (i_1,\dots,i_d) \in \Z^d.$ For $A \subset \R^d$ we let $\D_\delta(A)$ denote the intervals of $\D_\delta$ which intersect $A.$ Now let $X$ be a compactly supported random variable on $\R^d.$ Write the \textit{$\delta$-Shannon entropy} of $X$ by,  
\begin{equation}
    \Ed(X) = -\sum_{I \in \D_\delta}\pr(X\in I)\log \pr(X \in I),
\end{equation}
and the \textit{$\delta$-collision entropy} of $X$ by,
\begin{equation}
     \col(X) = -\log \sum_{I \in \D_\delta}\pr(X \in I)^2.
\end{equation}
We still have 
\begin{equation}\label{eq.mon}
\E_\delta(X) \geq \col(X),
\end{equation}
and the $\delta$-Shannon entropy inherits concavity from the finite setting. We also need the following facts.
\begin{lemma}[Continuity]\label{lem.contofent}
    Suppose that $X$ is a random variable on a compact set $A \subset \R^d$ and let $C > 1$, $ \delta > 0.$ Then
    $$\Ed(X) \leq \E_{C\delta}(X) +  O(1),$$
    with the implicit constant depending on $C$ and $d$ only.
    \begin{proof}
        Let $X_1$ be the random variable on $\D_\delta(A)$ which outputs the $I \in \D_\delta(A)$ for which $X \in I;$ let $X_2$ be the random variable on $\D_{C\delta}(A)$ which outputs the $J \in \D_{C\delta}(A)$ for which $X \in J.$ We have by the chain rule \eqref{eq.chainrule}, 
        $$\E(X_1,X_2) = \E(X_2) + \E(X_1 | X_2),$$
        which leads us to 
        $$\Ed(X) \leq \E_{C\delta}(X) + \E(X_1 | X_2).$$
        Finally, each $J \in \D_{C\delta}(\R^d)$ intersects at most $(\lceil C\rceil+2)^d$ cubes of $\D_\delta(\R^d)$. Hence, for each $J$ we have
        $$\E(X_1 | X_2 = J) \leq d\log (\lceil C\rceil+2),$$
        and so taking expectation gives us
        $$\E(X_1 |X_2) \leq d\log (\lceil C\rceil+2),$$
        and the result follows.
    \end{proof}
\end{lemma}
\begin{lemma}[Restriction]\label{lem.rest}
    Fix $\epsilon > 0$ and $\delta>0.$ Let $X$ be a random variable and suppose that $E$ is an event with $\pr(E) \geq 1-\epsilon.$ Then
    \begin{equation}
        (1-\epsilon)\Ed(X_{E}) \leq \Ed(X).
    \end{equation}
    Here $X_E$ is the random variable $X$ conditioned on $E.$
    \begin{proof}
        The distribution of $X$ is the mixture of the distributions of $X_E$ and $X_{E^c}$ with weights $\pr(E)$ and $\pr(E^c)$, respectively. By the concavity of Shannon entropy we therefore obtain
        \begin{equation}
        \pr(E)\Ed(X_E) + \pr(E^c)\Ed(X_{E^c}) \leq \Ed(X).
        \end{equation}
        A lower bound of the left-hand side is
        \begin{equation}
            (1-\epsilon)\Ed(X_E),
        \end{equation}
        and so we are done.
    \end{proof}
\end{lemma}
    \indent Let $C \geq 1$ and $s >0.$ We say that a random variable $X$ is \textit{$(s,C)$-Frostman} if $$\pr(X \in B(x,r)) \leq Cr^s$$ for all $x \in \R^d, r >0.$
\begin{lemma}
 Suppose that $X$ is $(s,C)$-Frostman. Then 
\begin{equation}
    \Ed(X) \geq s\log (1/\delta) - \log C - O(1).
\end{equation}
\begin{proof}
    We have
    \begin{align}
        \Ed(X) &= -\sum_{I \in \D_\delta}\pr(X \in I)\log \pr(X \in I)\\
        &\geq  -\sum_{I \in \D_\delta}\pr(X \in I)\log (C\delta^s) - O(1)\\
        &= s\log(1/\delta) - \log C - O(1).
    \end{align}
\end{proof}
\end{lemma}
We desire a submodular inequality in this setting too, which has previously been obtained in \cite{mator}. We include their proof for convenience.
\begin{lemma}[Discretised submodular inequality]\label{lem.dissubmod}
    Let $X,Y,Z,W$ be random variables taking values in  compact subsets of $\R^{k},\R^{l},\R^{m}, \R^n$ respectively. Fix $C > 1,$ $\delta > 0.$ Suppose each of the following:
    \begin{enumerate}
        \item If we know that the outcome of $X$ lies in $I \in \D_\delta(\R^k),$ then we are able to determine a choice of $2^m$ $J \in \D_{C\delta}(\R^m)$ which the outcome of $Z$ will lie;
        \item If we know that the outcome of $Y$ lies in $I \in \D_\delta(\R^l),$ then we are able to determine a choice of $2^m$ $J \in \D_{C\delta}(\R^m)$ which the outcome of $Z$ will lie;
        \item If we know the outcome of $X$ lies in $I \in \D_\delta(\R^k),$ and the outcome of $Y$ lies in $I' \in \D_\delta(\R^l),$ then we are able to determine a choice of $2^n$ $J \in \D_{C\delta}(\R^n)$ which the outcome of $W$ will lie.
        \end{enumerate}
        Then, $$\Ed(Z) + \Ed(W) \leq \Ed(X) + \Ed(Y) + O(1),$$
        where the implicit constant depends on $C,k,l,m,n$ only.
\begin{proof}
    Define the discrete random variables $X',Y'$ on the sample space $\D_\delta(\R^k), \D_\delta(\R^l)$ which output the $I \in \D_\delta(\R^k), J \in \D_\delta(\R^l)$ which the outputs of $X,Y$ lie in, respectively. Similarly, define the discrete random variables $Z',W'$ on the sample space $\D_{C\delta}(\R^m)$,$\D_{C\delta}(\R^n)$ respectively  which output the $I \in \D_{C\delta}(\R^m),J \in \D_{C\delta}(\R^n)$ which the outputs of $Z,W$ lie in, respectively. It is clear that
    $$\E(X') = \Ed(X), \qquad \E(Y') = \Ed(Y),$$
    and
    $$\E(Z') = \E_{C\delta}(Z), \qquad \E(W') = \E_{C\delta}(W).$$
    By construction,
    \begin{equation}
        \E(Z'|X')\leq m, \qquad \E(Z'|Y')\leq m,
        \qquad \E(W'|X',Y')\leq n.
    \end{equation}
    Apply Theorem \ref{thm.submod} with $Z'$ in the role of $X$, $(X',Y',Z')$ in the role of $Y$, and $(X',Z'),(Y',Z')$ in the roles of $Z,W$. This gives
    \begin{equation}
        \E(Z')+\E(X',Y',Z')
        \leq \E(X',Z')+\E(Y',Z').
    \end{equation}
    Using the chain rule and discarding the non-negative term $\E(Z'|X',Y')$ from the left-hand side, we obtain
    \begin{equation}
        \E(Z')+\E(X',Y')
        \leq \E(X')+\E(Y')+2m.
    \end{equation}
    Also,
    \begin{equation}
        \E(W')\leq \E(X',Y')+\E(W'|X',Y')
        \leq \E(X',Y')+n.
    \end{equation}
    Therefore
    $$\E(Z') + \E(W') \leq \E(X') + \E(Y') + 2m + n.$$
    Using the above identifications gives us,
    $$\E_{C\delta}(Z) + \E_{C\delta}(W) \leq \Ed(X) + \Ed(Y) +2m + n.$$
    Finally by the continuity of entropy (Lemma \ref{lem.contofent}) we have the result required.
\end{proof}
\end{lemma}
A useful rendition of this is the following.
\begin{lemma}\label{lem.sumprod}
    Let $X,Y,Z$ be i.i.d~random variables taking values in $[1,2].$ We have
    \begin{equation}
        \Ed((X+Y)Z) + 2\Ed(X) \leq \Ed(X+Y) + 2\Ed(XY) + O(1)
    \end{equation}
    for all $\delta >0.$ 
    \begin{proof}
   We essentially use the fact that the maps $(x,y) \mapsto xy, (x,y) \mapsto x+y$ are Lipschitz with constant $10$ (say), when restricted to $[1,4]^2,$ and the map $x \mapsto 1/x$ is Lipschitz with constant $10,$ when restricted to $[1,2].$
   
   \indent Using the facts above, suppose we know a ball of diameter $\delta$ in which $(X+Y,Z)$ lies. Then we are able to locate an interval of length $10\delta$ in which $(X+Y)Z$ will lie. Similarly, suppose we know a ball of diameter $\delta$ in which $(XZ,YZ)$ lies, then we are able to locate an interval of length $10\delta$ in which $(X+Y)Z$ will lie. 

   \indent Now suppose we jointly know two balls of diameter $\delta,$ where $(X+Y,Z)$ lies in one, and $(XZ,YZ)$ lies in the other. We therefore can determine an interval of length $\delta$ in which $Z$ lies, and so two intervals of length $10\delta,$ one containing $X$ and one containing $Y.$  

\indent Now, by Lemma \ref{lem.dissubmod}, we have
    \begin{equation}
        \Ed((X+Y)Z) + \Ed(X,Y,Z) \leq \Ed(X+Y,Z) + \Ed(XZ,YZ) + O(1). 
    \end{equation}
    Since $X,Y,Z$ are iid the result follows.
    \end{proof}
\end{lemma}
\section{Proof of Theorem 1.2}

\subsection{High multiplicity sets and a result of Orponen}

We recall some definitions motivated by \cite[Section 2.1]{orpahlfors}.

\begin{definition}
Let $K\subset\mathbb R^2$, $\theta\in S^1$, $N\geq 1$, and $\delta>0$.
Define the \emph{multiplicity function}
\[
    m_{K,\theta}(x,\delta)
    =
    N_\delta
    \big(
        K_\delta
        \cap
        \pi_\theta^{-1}(\pi_\theta(x))
    \big).
\]
Here $K_\delta$ is the $\delta$-neighbourhood of $K$ and $\pi_\theta$
is the orthogonal projection of $x$ to the line spanned by $\theta$.
Also write
\[
    H_\theta(K,N,\delta)
    =
    \{x\in\mathbb R^2:
    m_{K,\theta}(x,\delta)\geq N\}.
\]
\end{definition}

We will use the following consequence of the recent result of Orponen
\cite[Theorem 1.13]{orpahlfors}.

\begin{theorem}
For every $C,\epsilon,\sigma>0$, $s\in[0,1]$, and fixed $R\geq 1$,
there exists
\(
    \delta_0=\delta_0(C,\epsilon,\sigma,R)>0
\)
so that the following holds. Let $\mu$ be an $(s,C)$-regular measure
with
\(
    \operatorname{spt}\mu\subset B(0,R),
\)
and let $\nu$ be an $(\epsilon,C)$-Frostman measure on $S^1$. Then
\[
    \int_{S^1}
        \mu
        \big(
            H_\theta(\operatorname{spt}\mu,
            \delta^{-\sigma},\delta)
        \big)
    \,d\nu(\theta)
    \leq \epsilon
\]
for all $0<\delta<\delta_0$.
\end{theorem}

\begin{proof}
This follows from \cite[Theorem 1.13]{orpahlfors} by a fixed rescaling.
Indeed, let $a>0$, depending only on $R$, be sufficiently small so that
\[
    aB(0,R)\subset B(0,1/10),
\]
and apply the map $x\mapsto ax$. Write $\Delta=a\delta$.
Since the rescaled support is contained in $B(0,1/10)$, the
multiplicity restricted to a ball of radius $1$ appearing in
\cite[Theorem 1.13]{orpahlfors} agrees with the unrestricted multiplicity
above, up to harmless constant changes in scale.

Moreover,
\[
    \delta^{-\sigma}
    =
    a^\sigma\Delta^{-\sigma}
    \geq
    \Delta^{-\sigma/2}
\]
provided $\delta$ is sufficiently small. We may therefore apply
\cite[Theorem 1.13]{orpahlfors} with $\sigma/2$ in place of $\sigma$.
The rescaled measure remains regular with a constant depending only
on $C$ and $R$, and the result follows.
\end{proof}
\begin{remark}\label{rem.remark}
    To obtain Theorem \ref{thm.main} and Theorem \ref{thm.cor}, we would need the above for $s>1.$ In this case, $\sigma$ would be replaced by $\sigma + s-1.$ 
\end{remark}
We would like a result
where orthogonal projections are replaced with radial projections
with centres contained on a line. This is possible via a projective
transformation. We restate Definition 3.1 in this setting.

\begin{definition}
Let $K\subset\mathbb R^2$, $x\in\mathbb R^2$, $N\geq 1$, and
$\delta>0$. Define the \emph{radial multiplicity function}
\[
    \widetilde m_{K,x}(y,\delta)
    =
    N_\delta
    \big(
        K_\delta
        \cap
        \pi_x^{-1}(\pi_x(y))
    \big),
\]
where
\[
    \pi_x(y)=\frac{y-x}{|y-x|}
\]
is the radial projection to the circle of radius $1$ centred at $x$.
Also write
\[
    \widetilde H_x(K,N,\delta)
    =
    \{y\in\mathbb R^2\setminus\{x\}:
    \widetilde m_{K,x}(y,\delta)\geq N\}.
\]
\end{definition}

\begin{lemma}
For every $C,\epsilon,\sigma>0$ and $s\in[0,1]$, there exists
\(  \delta_0=\delta_0(C,\epsilon,\sigma)>0
\)
so that the following holds. Let $\mu$ be an $(s,C)$-regular measure
with
\(
    K:=\operatorname{spt}\mu\subset B(0,10),
\)
and let $\nu$ be an $(\epsilon,C)$-Frostman measure with
\(
    \operatorname{spt}\nu\subset l\cap B(0,10),
\)
where $l$ is a line. Suppose further that
\[
    1\leq
    \operatorname{dist}(K,\operatorname{spt}\nu)
    \leq 10,
    \qquad
    \operatorname{dist}(K,l)\geq \frac1{100}.
\]
Then
\[
    \int_l
        \mu
        \big(
            \widetilde H_x
            (K,\delta^{-\sigma},\delta)
        \big)
    \,d\nu(x)
    \leq \epsilon
\]
for all $0<\delta<\delta_0$.
\end{lemma}

\begin{proof}
We only apply this lemma when
\[
    l=\{0\}\times\mathbb R,
\]
so we complete the proof in this case. The general case follows by
a fixed Euclidean change of coordinates.

Define
\(
    P:\mathbb R^2\setminus l\longrightarrow\mathbb R^2
\)
by
\[
    P(x,y)=\left(\frac1x,\frac yx\right).
\]
Since
\[
    K\subset B(0,10),
    \qquad
    \operatorname{dist}(K,l)\geq\frac1{100},
\]
the map $P$ is bi-Lipschitz on a fixed neighbourhood of $K$.
Let $A\geq 1$ be a corresponding bi-Lipschitz constant.

For $t\in\mathbb R$ and $e=(e_1,e_2)\in S^1$ with $e_1\neq 0$,
let
\[
    l_t(e)=(0,t)+\operatorname{span}(e).
\]
For $r\neq 0$,
\[
\begin{aligned}
    P((0,t)+re)
    &=
    P(re_1,t+re_2)\\
    &=
    \left(\frac1{re_1},
    \frac{t+re_2}{re_1}\right)\\
    &=
    \left(0,\frac{e_2}{e_1}\right)
    +\frac1{re_1}(1,t).
\end{aligned}
\]
Therefore $P$ transforms lines through $(0,t)$ into lines parallel
to the vector $(1,t)$.

Define
\[
    Q(0,t)=\frac{(-t,1)}{\sqrt{1+t^2}}.
\]
Since
\[
    Q(0,t)^\perp=\operatorname{span}(1,t),
\]
the images under $P$ of lines through $(0,t)$ are fibres of the
orthogonal projection $\pi_{Q(0,t)}$.

The map $Q$ is bi-Lipschitz on $l\cap B(0,10)$. Consequently, if
\[
    \mu'=P_\#\mu,
    \qquad
    \nu'=Q_\#\nu,
\]
then $\mu'$ is $(s,C')$-regular and $\nu'$ is
$(\epsilon,C')$-Frostman, where $C'$ depends only on $C$ and
$\epsilon$. Moreover, $\operatorname{spt}\mu'$ is contained in a
fixed ball.

We now compare the multiplicities. Since $P$ is bi-Lipschitz on a
fixed neighbourhood of $K$, there is a constant $c_A>0$ such that
if
\[
    y\in K\cap
    \widetilde H_{(0,t)}
    (K,\delta^{-\sigma},\delta),
\]
then
\[
    P(y)\in
    H_{Q(0,t)}
    \big(
        P(K),
        c_A\delta^{-\sigma},
        A\delta
    \big).
\]
Indeed, the radial fibre through $y$ is mapped into a line parallel
to $(1,t)$, and covering numbers are preserved up to constants
depending only on $A$.

Write
\[
    \Delta=A\delta.
\]
For $\delta$ sufficiently small,
\[
    c_A\delta^{-\sigma}
    \geq
    \Delta^{-\sigma/2}.
\]
It follows that
\[
    P\big(
        K\cap
        \widetilde H_{(0,t)}
        (K,\delta^{-\sigma},\delta)
    \big)
    \subset
    H_{Q(0,t)}
    \big(
        P(K),
        \Delta^{-\sigma/2},
        \Delta
    \big).
\]
Therefore
\begin{equation*}
\begin{aligned}
    \int_l
        \mu
        \big(
            \widetilde H_x
            (K,\delta^{-\sigma},\delta)
        \big)
    \,d\nu(x)
    &\leq
    \int_{S^1}
        \mu'
        \big(
            H_\theta
            (P(K),\Delta^{-\sigma/2},\Delta)
        \big)
    \,d\nu'(\theta).
\end{aligned}
\end{equation*}
Applying Theorem 3.2, with $\sigma/2$ in place of $\sigma$, gives
the result.
\end{proof}

We prove an entropic version of this below. Fix $z\in\mathbb R$.
We wish to compare the radial projection $\pi_{(0,z)}$ with the map
\[
    \angle_z:\mathbb R^2\setminus\{x = 0\}
    \longrightarrow\mathbb R,
    \qquad
    \angle_z(x,y)=\frac{y-z}{x},
\]
which measures the slope between $(x,y)$ and $(0,z)$.

If $|x|\geq 1/2$ and $x,y,z\in[-2,2]$, then
$|\angle_z(x,y)|\leq 8$, and
\begin{equation}\label{eq.slope-radial}
    \pi_{(0,z)}(x,y)
    =
    \operatorname{sgn}(x)
    \frac{(1,\angle_z(x,y))}
    {\sqrt{1+\angle_z(x,y)^2}}.
\end{equation}
Consequently, on each of the half-planes $x>0$ and $x<0$, the slope
and radial projection parametrisations are bi-Lipschitz comparable,
with absolute constants.

\begin{proposition}
Let $C\geq 1$, $0<s\leq 1/2$, $0<\epsilon<s$, and $\sigma>0$.
There exists
\(
    \delta_0=\delta_0(C,\epsilon,\sigma)>0
\)
so that the following holds. Let $\mu,\nu$ be $(s,C)$-regular
measures on $[-2,2]$ and let $\xi$ be an $(s,C)$-Frostman measure
on $[-2,2]$. Suppose that
\[
    \operatorname{dist}
    \big(
        \operatorname{spt}\mu\times\operatorname{spt}\nu,
        \{0\}\times\operatorname{spt}\xi
    \big)
    \geq 1,
\]
and
\[
    \operatorname{dist}(\operatorname{spt}\mu,0)\geq\frac12.
\]
Let $X,Y,Z$ be independent random variables distributed by
$\mu,\nu,\xi$ respectively. Then
\[
    H_\delta
    \left(
        \frac{Y-Z}{X}
    \right)
    \geq
    2(1-\epsilon)(s-\sigma)\log(1/\delta)
    -O_C(1).
\]
\end{proposition}

\begin{proof}
Fix $0<\delta<\delta_0$. Since $X,Y$ are independent of $Z$, the
distribution of $(Y-Z)/X$ is the mixture, with respect to $\xi$, of
the distributions of $(Y-z)/X$. By the concavity of Shannon entropy,
applied to the fixed partition $\mathcal D_\delta(\mathbb R)$,
\[
    H_\delta\left(\frac{Y-Z}{X}\right)
    \geq
    \int H_\delta\left(\frac{Y-z}{X}\right)\,d\xi(z).
\]

Fix $z\in\operatorname{spt}\xi$. Write
\[
    K=\operatorname{spt}\mu\times\operatorname{spt}\nu
\]
and
\[
    M_z
    =
    (\mu\times\nu)
    \big(
        \widetilde H_{(0,z)}
        (K,\delta^{-\sigma},\delta)
    \big).
\]
Also write
\[
    E_z
    =
    K\setminus
    \widetilde H_{(0,z)}
    (K,\delta^{-\sigma},\delta),
    \qquad
    a_z=1-M_z.
\]
If $a_z=0$, the lower bound we obtain below is trivial. We therefore
assume that $a_z>0$ and let $\rho_z$ be the probability measure
\[
    \rho_z
    =
    \frac{1}{a_z}
    (\mu\times\nu)|_{E_z}.
\]
Let $(X',Y')$ be distributed according to $\rho_z$.

We first lower-bound the collision entropy of
\[
    \frac{Y'-z}{X'}.
\]
Fix a sufficiently small absolute constant $c_0>0$, and let
$\mathcal D_{c_0\delta}(\mathbb R)$ be the standard partition of
$\mathbb R$ into intervals of length $c_0\delta$. For
$I\in\mathcal D_{c_0\delta}(\mathbb R)$, write
\[
    T_I
    =
    \{(x,y)\in K:\angle_z(x,y)\in I\}.
\]
By \eqref{eq.slope-radial}, these are unions of at most two tubes of width comparable
to $\delta$, one in each of the half-planes $x>0$ and $x<0$.

We claim that
\[
    E_z\cap T_I
\]
can be covered by $\lesssim\delta^{-\sigma}$ balls of radius
comparable to $\delta$. To see this, consider first the part with
$x>0$, and suppose it is non-empty. Choose
\[
    p=(x_0,y_0)\in E_z\cap T_I
\]
with $x_0>0$, and write
\[
    m_0=\angle_z(p).
\]
If $q=(x,y)\in E_z\cap T_I$ also has $x>0$, then
\[
    |\angle_z(q)-m_0|
    \leq c_0\delta.
\]
Hence the point
\[
    q'=(x,z+xm_0)
\]
lies on the same radial fibre through $(0,z)$ as $p$, and
\[
    |q-q'|
    \leq 2c_0\delta.
\]
Choosing $c_0$ sufficiently small gives $q'\in K_\delta$.
Since $p\notin
\widetilde H_{(0,z)}(K,\delta^{-\sigma},\delta)$,
the set
\[
    K_\delta
    \cap
    \pi_{(0,z)}^{-1}(\pi_{(0,z)}(p))
\]
can be covered by fewer than $\delta^{-\sigma}$ balls of radius
$\delta$. It follows that the part of $E_z\cap T_I$ with $x>0$
can be covered by $\lesssim\delta^{-\sigma}$ such balls.
The same argument applies in the half-plane $x<0$, proving the claim.

Since $\mu\times\nu$ is $(2s,C')$-Frostman, every ball of radius
comparable to $\delta$ has $\mu\times\nu$-measure
$\lesssim_C\delta^{2s}$. Therefore
\[
    \rho_z(T_I)
    \lesssim_C
    a_z^{-1}\delta^{2s-\sigma}.
\]
Also, since $K$ is upper $(2s,C')$-regular and the map
$\angle_z$ is Lipschitz on the region under consideration, the
number of intervals $I$ for which $\rho_z(T_I)>0$ is
\[
    \lesssim_C\delta^{-2s}.
\]
Consequently,
\[
\begin{aligned}
    \sum_I\rho_z(T_I)^2
    &\lesssim_C
    \delta^{-2s}
    a_z^{-2}\delta^{4s-2\sigma}\\
    &=
    a_z^{-2}\delta^{2s-2\sigma}O_C(1).
\end{aligned}
\]
It follows that
\[
    \operatorname{col}_{c_0\delta}
    \left(
        \frac{Y'-z}{X'}
    \right)
    \geq
    2(s-\sigma)\log(1/\delta)
    +2\log a_z
    -O_C(1).
\]

By the restriction estimate, the monotonicity of entropy, and
Lemma 2.4,
\[
\begin{aligned}
    H_\delta
    \left(
        \frac{Y-z}{X}
    \right)
    &\geq
    H_{c_0\delta}
    \left(
        \frac{Y-z}{X}
    \right)-O(1)\\
    &\geq
    a_z
    H_{c_0\delta}
    \left(
        \frac{Y'-z}{X'}
    \right)-O(1)\\
    &\geq
    a_z
    \operatorname{col}_{c_0\delta}
    \left(
        \frac{Y'-z}{X'}
    \right)-O(1)\\
    &\geq
    2a_z(s-\sigma)\log(1/\delta)
    +2a_z\log a_z
    -O_C(1).
\end{aligned}
\]
Since
\[
    x\log x\geq -O(1),
    \qquad 0\leq x\leq 1,
\]
we obtain
\[
    H_\delta
    \left(
        \frac{Y-z}{X}
    \right)
    \geq
    2(1-M_z)(s-\sigma)\log(1/\delta)
    -O_C(1).
\]

Finally, $\mu\times\nu$ is $(2s,C')$-regular and $2s\leq 1$.
Since $\epsilon<s$, the measure $\xi$ is also
$(\epsilon,C')$-Frostman. Lemma 3.4 therefore gives
\[
    \int M_z\,d\xi(z)\leq\epsilon.
\]
Integrating the preceding entropy estimate gives
\[
    \int H_\delta\left(\frac{Y-z}{X}\right)\,d\xi(z)
    \geq
    2(1-\epsilon)(s-\sigma)\log(1/\delta)
    -O_C(1),
\]
and the result now follows from the concavity estimate at the start
of the proof.
\end{proof}

\begin{corollary}
Let $\mu$ be an $(s,C)$-regular measure on $[1,2]$, where
$0<s\leq 1/2$. Let $X,Y,Z$ be iid random variables distributed
by $\mu$. For every $0<\epsilon<s$ and $\sigma>0$, there exists
$\delta_0>0$ such that
\[
    H_\delta((X+Y)Z)
    \geq
    2(1-\epsilon)(s-\sigma)\log(1/\delta)
    -O_C(1)
\]
for all $0<\delta<\delta_0$.
\end{corollary}

\begin{proof}
Let
\[
    \mu_1=(x\mapsto x^{-1})_\#\mu,
    \qquad
    \nu_1=\mu,
    \qquad
    \xi_1=(x\mapsto -x)_\#\mu.
\]
The map $x\mapsto x^{-1}$ is bi-Lipschitz on $[1,2]$, so these
measures satisfy the hypotheses of Proposition 3.6 after increasing
the regularity constant by an absolute factor.

Let $X_1,Y_1,Z_1$ be independent random variables distributed by
$\mu_1,\nu_1,\xi_1$, respectively. Proposition 3.6 gives
\[
    H_\delta\left(\frac{Y_1-Z_1}{X_1}\right)
    \geq
    2(1-\epsilon)(s-\sigma)\log(1/\delta)
    -O_C(1).
\]
If $A,B,C$ are iid  with distribution $\mu$, then
\( (Y_1-Z_1)/X_1 \) has the same distribution as \(A(B+C)\). By
symmetry this has the same distribution as \(C(A+B)\).
The result follows.
\end{proof}

\subsection{Proof of Theorem 1.2}

\begin{proof}[Proof of Theorem 1.2]
Fix $0<s\leq 1/2$ and $\eta>0$. Choose
$0<\epsilon<s$ and $\sigma>0$ sufficiently small that
\[
    2(1-\epsilon)(s-\sigma)\geq 2s-\eta.
\]
Let $\delta_0$ be that given by Corollary 3.7. By Corollary 3.7,
\[
    H_\delta((X+Y)Z)
    \geq
    (2s-\eta)\log(1/\delta)-O_C(1).
\]
By Lemma 2.6,
\[
    H_\delta(X)
    \geq
    s\log(1/\delta)-\log C-O(1).
\]
Therefore
\[
\begin{aligned}
    (4s-\eta)\log(1/\delta)-O_C(1)
    &\leq
    H_\delta((X+Y)Z)+2H_\delta(X)\\
    &\leq
    H_\delta(X+Y)+2H_\delta(XY)+O_C(1),
\end{aligned}
\]
where the second inequality follows from Lemma 2.8. This proves the
result.
\end{proof}
\section{Sharpness of Theorem \ref{thm.main}}\label{section.sharp}
We use the language of iterated function systems, see \cite[Chapter 11]{falbk}. Fix $0 < s \leq 1/2, \eta >0.$ Let $P \subset [0,1)$ be an arithmetic progression of length $N,$ starting from $0,$ and of step $1/N.$ Consider the iterated function system $\calF = \{cx + p\}_{p \in P},$ where $c \leq 1/N$ is chosen so that $s = \tfrac{\log N}{\log 1/c}.$ Let $A$ be the attractor of $\calF$ and let $\mu$ be the self-similar measure with uniform weights on $P.$ It is well known that $\mu$ is $(s,C)$-regular for some $C = C(s,N) >0.$ Since $P+P$ has cardinality $2N-1,$ the set $A+A$ is the attractor of the iterated function system with contraction ratio $c$ and translations in $P+P.$ Hence
\begin{equation}
    \overline{\dim}_{\rm B}(A+A)
    \leq \frac{\log(2N-1)}{\log(1/c)}
    =s\frac{\log(2N-1)}{\log N}.
\end{equation}
Choose $N$ large enough that the right-hand side is less than $s+\eta/8.$ It follows that, for all $\delta>0$ small enough,
\begin{equation}
    N_\delta(A+A)\leq \delta^{-s-\eta/4}.
\end{equation}
Now $C$ will depend on $s$ and $\eta.$ Consider the map $x \rightarrow 2^x.$ This restricted to $[0,1]$ is bi-Lipschitz and its image is contained in $[1,2].$ The image measure $\nu$ of $\mu$ under the map is therefore still $(s,C')$-regular, where $C'$ depends on $C$ only. If $A'=\spt\nu$, then $A'A'=2^{A+A}.$ Since $x\mapsto 2^x$ is bi-Lipschitz on $[0,2]$, after absorbing the fixed change of scale, for all $\delta>0$ small enough we have
\begin{equation}
    N_\delta(A'A')\leq \delta^{-s-\eta/2}.
\end{equation}
Also, the upper regularity of $A'$ gives $N_\delta(A')\lesssim_{C'}\delta^{-s}$, and hence
\begin{equation}
    N_\delta(A'+A')\lesssim_{C'}\delta^{-2s}.
\end{equation}
Let $X,Y$ be iid random variables distributed by $\nu.$ For such a $\delta$ we have
\begin{align*}
    \Ed(X+Y) + 2\Ed(XY)
    &\leq \log N_\delta(A'+A')N_\delta(A'A')^2+O(1)\\
    &\leq (4s + \eta)\log (1/\delta)+O_C(1). 
\end{align*}
Thus Theorem \ref{thm.main} and \eqref{eq.sharp} are both sharp. 
\bibliographystyle{alpha}
\bibliography{references.bib}
\end{document}